\documentclass[aps,prl,twocolumn,groupedaddress]{revtex4-1}

\usepackage{amsmath}           
\usepackage{graphicx}          
\usepackage{kotex}             
\usepackage{dcolumn}           
\usepackage{here}
\usepackage{amssymb}
\usepackage{centernot}
\usepackage{mathtools}
\usepackage{stmaryrd}

\newcommand{\xMapsto}[2][]{\ext@arrow 0599{\Mapstofill@}{#1}{#2}}
\def\Mapstofill@{\arrowfill@{\Mapstochar\Relbar}\Relbar\Rightarrow}

\usepackage[utf8]{inputenc}
\usepackage{listings}
\usepackage{lipsum}
\usepackage{subfig}

\newtheorem{theorem}{Theorem}[section]

\newenvironment{Even Case}[1][Even Case]{\begin{trivlist}
		\item[\hskip \labelsep {\bfseries #1}]}{\end{trivlist}}
\newenvironment{Odd Case}[1][Odd Case]{\begin{trivlist}
		\item[\hskip \labelsep {\bfseries #1}]}{\end{trivlist}}

\newcommand{\qed}{\nobreak \ifvmode \relax \else
	\ifdim\lastskip<1.5em \hskip-\lastskip
	\hskip1.5em plus0em minus0.5em \fi \nobreak
	\vrule height0.75em width0.5em depth0.25em\fi}

\usepackage{mathtools}
\newcommand\Set[2]{\{\,#1\mid#2\,\}}

\begin{document}               

\title{Numerical Approach for Fermat's last theorem}        

\author{Youngik Lee}            

\affiliation{Dept. of Physics, University of Bonn, Bonn-Cologne Graduate School, Germany}



\begin{abstract}
    This research focuses on the Numerical approach for Fermat's Last theorem. We can induce an Alternative form of Fermat's last theorem by using particular geometric mapping $\mathcal{M}$ on a Cartesian plane to a Torus. It transforms the Fermat's Last Theorem to finding a rational cross point between two parametric curves on the torus. In the end, this research shows the movement of the point, on the line $x^n+y^n=1$, has an acceleration phase transition near ($x,n$)=(0,2). Moreover, the studies about the relationship between this acceleration transition and the solution for the Fermat's Diophantine equation in the case of $n>$2, need further investigation.

\end{abstract}

\maketitle			        
\tableofcontents

\section{1. Introduction}

The Fermat's last theorem is a famous problem in number theory, and proved was formally published by Andrew Wiles in 1995 by showing the Taniyama–Shimura conjecture is valid, which states that elliptic curves over the field of rational numbers are related to modular forms.
The problem was suggested by Pierre de Fermat in 1637, and famous with the note that he had a proof that was too large to fit in the margin.

Fermat's last theorem is simply,

\begin{theorem}
	\emph{Fermat's Last theorem}
	
	The Diophantine equation,
	\begin{equation}\label{eq7}
	x^n+y^n=z^n
	\end{equation}
	has no integer solutions for n$>$2 and x,y,z $\neq$ 0. \qed
\end{theorem}

Here, the Diophantine equation is an equation in which only integer solutions are allowed.

Furthermore, this paper suggests a few Alternative forms of Fermat's last theorem and show the theorem is equivalent to finding the rational cross point between two parametrized curves on the torus, by using the particular mapping from Cartesian coordinate. Moreover, we can define the concept of acceleration of the point exists on the line, $x^n+y^n=1$, and can show the point acceleration has a phase transition near $n$=2, $x$=0.

\section{2. Theory}

\subsection{2.1 Coordinate Rescaling}

When we divide out $z^n$ in both side of Eq. (1), we can change the Eq. (1) into Eq. (2),

\begin{equation}\label{eq7}
\left(\frac{x}{z}\right)^n+\left(\frac{y}{z}\right)^n=1
\end{equation}

Then we can rewrite \textbf{Theorem .1} in Alternative form.

\begin{theorem}
	\emph{Alternative Fermat's Last theorem I}
	
	When X, Y$\in \mathbb{Q}$, X,Y$\in [0,1]$, the equation, 
	\begin{equation}\label{eq7}
	X^n+Y^n=1
	\end{equation}
	has no rational solutions $^{\forall}n>2$ \qed
\end{theorem}

Fig. 1 shows the set of graphs for Eq. (3) when $n$ changes from 1 to 10.
We can see that when $n$ goes larger, the curves are more converge to the unit square. Now let us call this curve as Fermat's Curve.

\begin{figure}[H]
	\centering
	\includegraphics[angle=0, width=6cm, height=6cm]{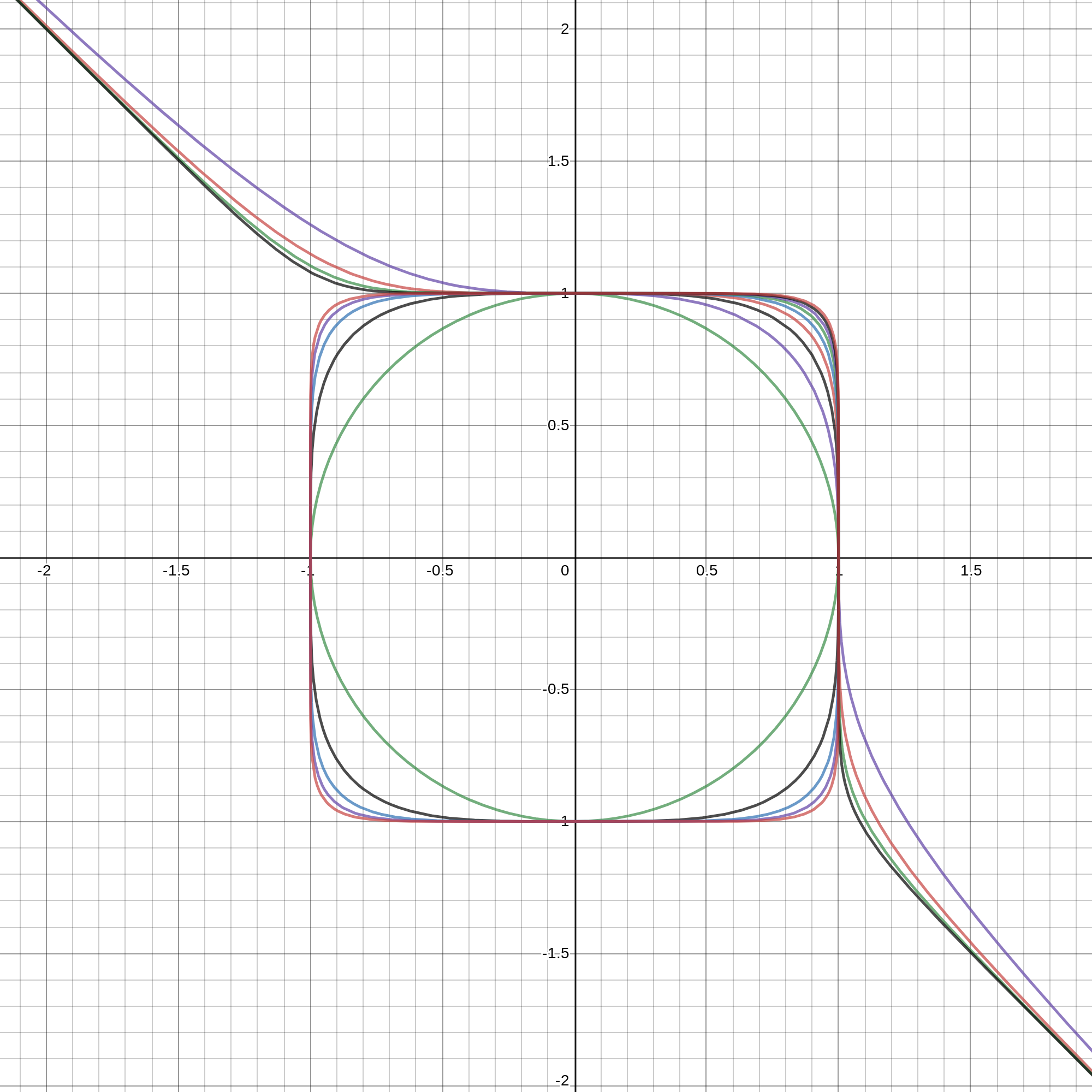}
	\caption{\label{fg3} The set of Fermat's Curve, with different $n$ from 1 to 10. The figure shows when $n$ increase, the curves converge to unit square.}
\end{figure}

Now we can say that alternative Fermat's theorem is equivalent to proving the existence of any rational cross point between two set $\mathcal{F}$ and $\mathcal{D}$, which are defined in \textbf{Theorem .3}.

(Here, a rational point means the point which has the only rational number as its coordinate.)

\begin{theorem}
	\emph{Alternative Fermat's Last Theorem II}
	
	When the set $\mathcal{F}$ and $\mathcal{D}$ are defined in the range of $^\forall$x,y $\in$ [0,1] like below,
	
	\begin{align*}
	\mathcal{F}  &\coloneqq \Set{(x, y)}{\operatorname{}x^n+y^n=1, ^\forall x, y\in \mathbb{Q}}\\
	\mathcal{D}  &\coloneqq \Set{(x, y)}{\operatorname{}x\times y, ^\forall x, y\in \mathbb{Q}}
	\end{align*}
	
	Then the Fermat's Last Theorem is equivalent to the below equation.
	
	\begin{equation}\label{eq7}
	\mathcal{F}\cap\mathcal{D}=\phi, ^{\forall}n>2 \qed
	\end{equation}
\end{theorem}

\subsection{2.2 Set of lines and positive integer plane, $\mathcal{U}$}

First, let us consider another expression of $\mathcal{D}$ by using linear mapping. We can define new set $\mathbb{N'}$ like Eq. (5)

\begin{equation}\label{eq7}
\mathbb{N'}=\mathbb{N} \cup \{0\}
\end{equation}

Then we can express the positive integer plane $\mathcal{U}$, like Eq. (6)

\begin{align}
\mathcal{U}  &\coloneqq \Set{(x, y)}{\operatorname{}x\times y, ^\forall x, y\in \mathbb{N'}}
\end{align}

There are various ways to express $\mathcal{U}$ as a set of different geometric objects, but in this case, let us use the line as an object element.

\begin{align}
l(a, b, n)  &\coloneqq \Set{(x, y)}{x=an, y=bn, ^\forall a,b\in \mathbb{N'}}
\end{align}

Here, $n \in \mathbb{N}$. Then it is straight forward to show, the whole union of $l$ satisfies Eq. (8)

\begin{equation}
\cup ^\forall l  \equiv \mathcal{U}
\end{equation}

Here $l$ is discrete line which have period of $n$, with unit length $R=\sqrt{a^2+b^2}$

Therefore as we can also see in Fig. 2, \textbf{Theorem .3} is equivalent to asking whether is there a rational cross point between line $l$ and curve set $\mathcal{F}$.

\begin{theorem}
	\emph{Alternative Fermat's Last Theorem III}
	
	In the range of $^\forall$x,y $\in$ [0,1], when set $\mathcal{F}$ and line l is defined  like below,
	
	\begin{align*}
	\mathcal{F}  &\coloneqq \Set{(x, y)}{\operatorname{}x^n+y^n=1, ^\forall x, y\in \mathbb{Q}}\\
	l'(a, b, n)  &\coloneqq \Set{(x, y)}{x=an, y=bn, ^\forall a, b\in \mathbb{Q}_{[0,1]}}
	\end{align*}
	
	Here, $n \in \mathbb{N}$.
	Then the Fermat's Last Theorem is equivalent to the below equation.
	
	\begin{equation}\label{eq7}
	\mathcal{F}\cap l'=\phi, ^{\forall}n>2 \qed
	\end{equation}
\end{theorem}

To prove the \textbf{Theorem .4}, the next step is finding the geometrical object $\mathcal{T}$, which we can project the line $l'$, and also conserve its periodicity by $n$.

\begin{figure}[H]
	\centering
	\includegraphics[angle=0, width=6cm, height=6cm]{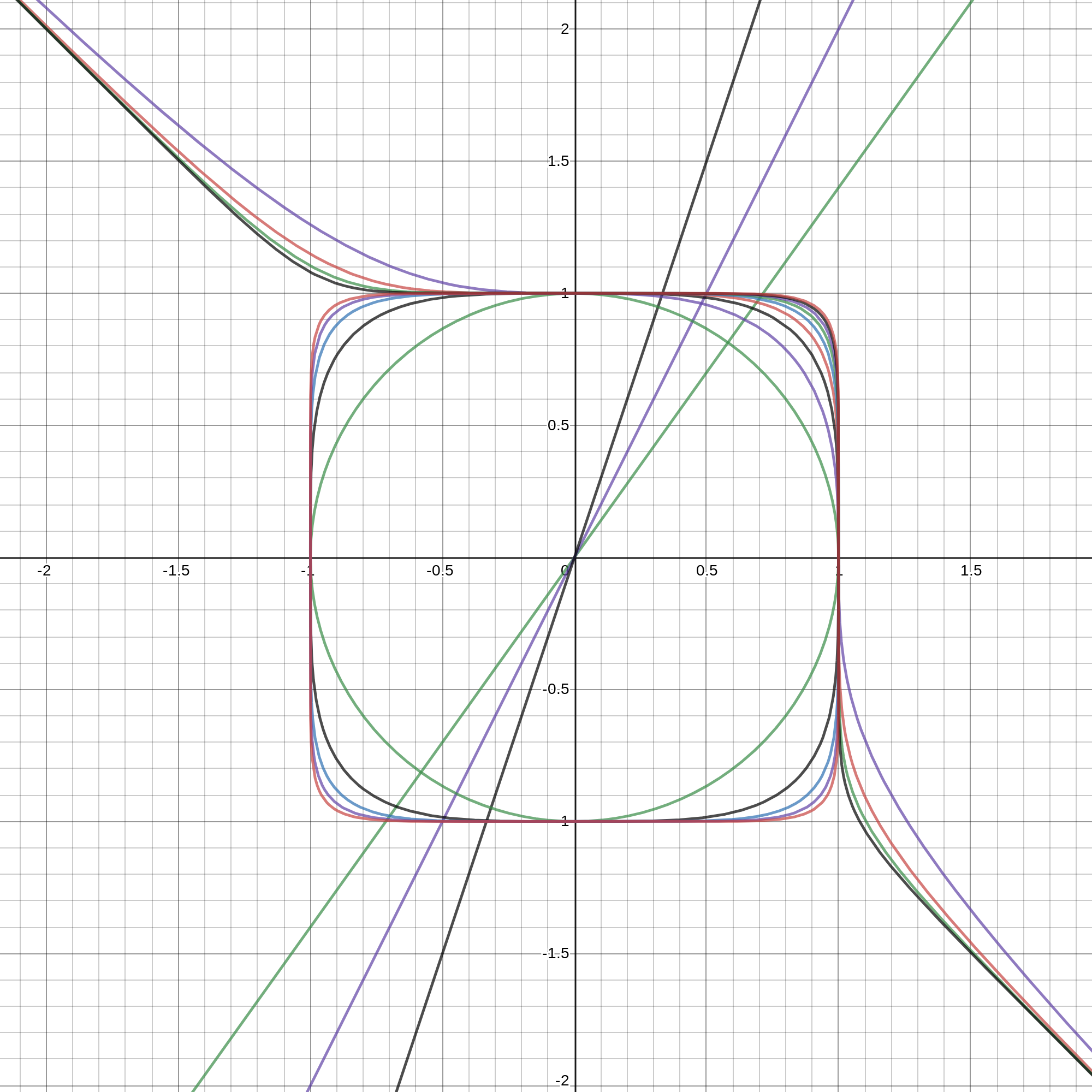}
	\caption{\label{fg3} The illustration for $\mathcal{F}$  and lines on plane. According to the \textbf{Theorem .4}, Fermat's last theorem is equivalent to find the rational crossing point between straight line and curves on the figure.}
\end{figure}

\subsection{2.3 Mapping on the Torus}

From now, the geometric object we are going to study is Torus. The reason for this is, Torus is a nice geometrical object which can preserve the "periodicity" after mapping other linear function on it.
The equation of Torus is defined like Eq. (10)

\begin{align}
\mathcal{T}  &\coloneqq \Set{(x, y)}{\operatorname{}(\sqrt{x^2+y^2}-R)^2+z^2=r^2 }
\end{align}

Moreover, we can write the general coordinate of Torus as Eq. (11)

\begin{equation}\label{eq7}
\sigma(u,v)=
\begin{cases}
(R+r\cos v)\cos u& \text{ x-corrdinate }\\ 
& \\
(R+r\cos v)\sin u& \text{ y-corrdinate } \\
& \\
r\sin v& \text{ z-corrdinate }
\end{cases}
\end{equation}

Here let us assume R$>$r, and Fig. 3 shows the direction of $u$ and $v$ coordinate on Torus.

\begin{figure}[H]
	\centering
	\includegraphics[angle=0, width=6.5cm, height=3.6cm]{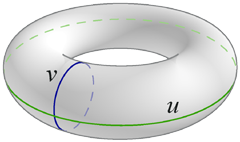}
	\caption{\label{fg3} The two degree of freedom for angular movement on Torus, u and v.}
\end{figure}

To get the parametric  curve equation on the torus, first set initial starting point at $(x_0,y_0)$ and move in the direction of $a\sigma_u+b\sigma_v$. The parametric equation we are going to use is Eq. (12).

\begin{equation}\label{eq7}
\begin{aligned}
u=at+x_0  \\
v=bt+y_0
\end{aligned}
\end{equation}

Thus, the parametric equation of the curve on the Torus becomes Eq. (13).

\begin{equation}\label{eq7}
\gamma(t)=\sigma(u(t),v(t))
\end{equation}

Therefore, when we set $(x_0,y_0)$=(0,0), then $\gamma(t)$ is simply given as $\sigma(at,bt)$.

If $b\neq$0 and $a$, $b$ is rational, the parameterization above is periodic, and the curve is closed. 
However, in a case when $a$, $b$ irrational, the curve is not closed, and it becomes a dense subset of the Torus. (If $b$=0, the curve is closed.)

When the parametrized curve is closed, we can define the period "T," which satisfy Eq. (14)

\begin{equation}\label{eq7}
\begin{aligned}
u(T)=u(0)  \\
v(T)=v(0)
\end{aligned}
\end{equation}

Furthermore, by using T, we can define disconnected parametrized cure $\gamma'$ as Eq. (15)

\begin{align}
\gamma'  &\coloneqq \Set{(x, y, n)}{(\mathbb{F}(\gamma(nT)),n),  n\in\mathbb{N})}
\end{align}

Here, $\mathbb{F}(\vec{x})$, and $proj_i(\vec{x})$ are defined as Eq. (16)

\begin{equation}\label{eq7}
\begin{cases}
\mathbb{F}(\vec{x}) = proj_1(\vec{x})\times proj_2(\vec{x})& \text{}\\ 
& \\
proj_i(\vec{x})=x_i& \text{ }
\end{cases}
\end{equation}

And with Eq. (13), we can rewrite Eq. (15) as,

\begin{align}
\gamma'  &\coloneqq \Set{(x, y, n)}{(\mathbb{F}(\sigma(u(nT),v(nT)),n),  n\in\mathbb{N})}
\end{align}

Thus finally we can get a simple form of $\gamma'$

\begin{align}
\gamma'  &\coloneqq \Set{(x, y, n)}{((R+r),0,n)}
\end{align}

And it induces one-to-one mapping $\mathcal{M}$, between cartesian and torus coordinate, like Eq. (19)

\begin{align}
\mathcal{M}  &\coloneqq l(a,b,n) \mapsto \gamma'(R, r, n)
\end{align}

Now let us visualize the line mapping process on a Torus. For example in $b/a$=3 case, Fig. 4 shows the process of infinite line mapping into a unit square, start with the top left illustration, in a clockwise direction.

\begin{figure}[H]
	\centering
	\includegraphics[angle=0, width=6cm, height=6cm]{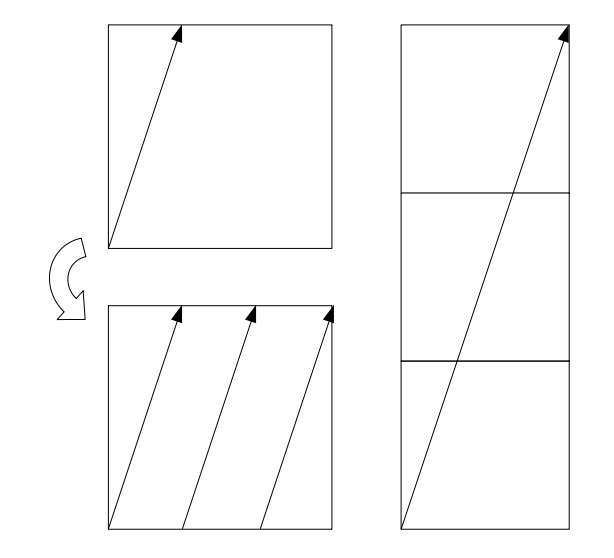}
	\caption{\label{fg3} The process of projecting infinite line on unit square. Start with the top left, and proceed in clockwise direction.}
\end{figure}

When $b/a$ gives integer value, infinitely extrapolated line always cross the point $(b,a)$. Then we can overlap all line paths on $a\times b$ squares in a unit square like the bottom left illustration of Fig. 4.
After that if you rolled up this unit square like Fig. 5, you can get the open end cylinder.

\begin{figure}[H]
	\centering
	\includegraphics[angle=0, width=5.5cm, height=3cm]{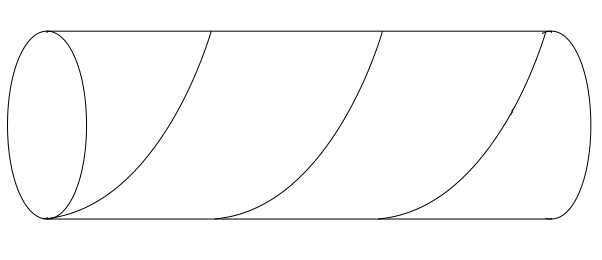}
	\caption{\label{fg3} The open end cylinder with curved lines.}
\end{figure}

At last, when you connect the two open end of the cylinder, one can make a Torus, like Fig. 6.

\begin{figure}[H]
	\centering
	\includegraphics[angle=0, width=6.5cm, height=3.6cm]{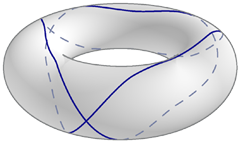}
	\caption{\label{fg3} The line mapping on Torus. $b/a$=3 case.}
\end{figure}

This is the simple visualization of the mapping $\mathcal{M}$

Then we can ask that whether is it available to describe $\gamma(t)$ on torus in explicit form. To answer this question, let us take a look at the explicit equation form of a mapped line on Torus by calculating geodesics.

\subsection{2.4 The Geodesic Equation on Torus}

To calculate the period T, we need to define the velocity and acceleration of u, v. Take the partial derivatives of this parameterization, then we can get Eq. (20), Eq. (21)

\begin{equation}\label{eq7}
\vec{x}_u=
\begin{cases}
(R+r\cos v)\sin u& \text{ x-corrdinate }\\ 
& \\
(R+r\cos v)\cos u& \text{ y-corrdinate } \\
& \\
0& \text{ z-corrdinate }
\end{cases}
\end{equation}

\begin{equation}\label{eq7}
\vec{x}_v=
\begin{cases}
-r\cos u\sin v& \text{ x-corrdinate }\\ 
& \\
-r\sin u\sin v& \text{ y-corrdinate } \\
& \\
r\cos v& \text{ z-corrdinate }
\end{cases}
\end{equation}

Then we can compute inner products to find the coefficients of the first fundamental form, $E, G, F$ like Eq. (22).

\begin{equation}\label{eq7}
\begin{aligned}
E = \vec{x}_u \cdot \vec{x}_u = (R + a \cos v)^2  \\
F = \vec{x}_u \cdot \vec{x}_v = 0 \\
G = \vec{x}_v \cdot \vec{x}_v=a^2
\end{aligned}
\end{equation}

This gives us the line element $ds^2 = (R + a \cos v)^2 du^2 + r^2 dv^2$ , from which we can get the metric, $g_{ij}$

\begin{equation}\label{eq7}
g_{ij} = \begin{bmatrix}
(R + r \cos v)^2 & 0\\ 
0 & r^2
\end{bmatrix}
\end{equation}

\begin{equation}\label{eq7}
g_{ij, u} = \begin{bmatrix}
0 & 0\\ 
0 & 0
\end{bmatrix}
\end{equation}

\begin{equation}\label{eq7}
g_{ij, v} = \begin{bmatrix}
-2r\sin v(R + r \cos v) & 0\\ 
0 & 0
\end{bmatrix}
\end{equation}

To calculate the geodesic equation, first we need to know the value for the nonzero Christoffel symbols of the second kind.

Christoffel symbols $\Gamma^i_{kl}$ is defined like, Eq. (26)

\begin{equation}\label{eq7}
\begin{aligned}
\Gamma^i_{kl} = \frac{1}{2}g^{im}(g_{mk,l}+g_{ml,k}-g_{kl,m})
\end{aligned}
\end{equation}

And for $u$ and $v$ coordinate, the non-zero value for Christoffel symbols of the second kind is $\Gamma^u_{uv}$ and $\Gamma^v_{uu}$.

\begin{equation}\label{eq7}
\begin{aligned}
\Gamma^u_{uv} = \frac{1}{2}[g^{uu}(g_{uv,u}+g_{uu,v}-g_{uu,u})  \\
+g^{uv}(g_{vu,u}+g_{vu,u}-g_{uu,v})]
\end{aligned}
\end{equation}

\begin{equation}\label{eq7}
\begin{aligned}
\Gamma^v_{uu} = \frac{1}{2}[g^{vu}(g_{vu,u}+g_{vu,u}-g_{uu,v})  \\
+g^{vv}(g_{vu,u}+g_{vu,u}-g_{uu,v})]
\end{aligned}
\end{equation}

By using Eq. (23)$\sim$(25), Eq. (27)$\sim$(28) can be change to,

\begin{equation}\label{eq7}
\Gamma^u_{uv} = -\frac{r\sin v}{(R+a\cos v)^2}
\end{equation}

\begin{equation}\label{eq7}
\Gamma^v_{uu} = \frac{1}{a}\sin v(R+a\cos v)
\end{equation}

After calculate the non-zero Christoffel values, we can solve geodesics equation Eq. (31)

\begin{equation}\label{eq7}
\ddot{x}^a +\Gamma_{bc}^a \dot{x}^b \dot{x}^c = 0
\end{equation}

From Eq. (31), we can directly get, 

\begin{equation}\label{eq7}
\begin{cases}
\ddot{u}+2\Gamma^{u}_{uv} \dot{u}\dot{v} = 0& \text{ }\\ 
& \\
\ddot{v}+2\Gamma^{v}_{uu} \dot{u}^2 = 0& \text{ }
\end{cases}
\end{equation}

Moreover, by using the non-zero Christoffel values(Eq. (29)$\sim$(30)) we can get Eq. (33)

\begin{equation}\label{eq7}
\begin{cases}
\ddot{u}-\frac{2r\sin v}{(R+a\cos v)^2}\dot{u}\dot{v} = 0& \text{ }\\ 
& \\
\ddot{v}+\frac{1}{r} \sin v (R + r \cos v) \dot{u}^2 =0& \text{ }
\end{cases}
\end{equation}

\begin{figure*}[htp]
	\centering
	\subfloat[ velocity, n=1.9]{\label{figur:2}\includegraphics[width=50mm]{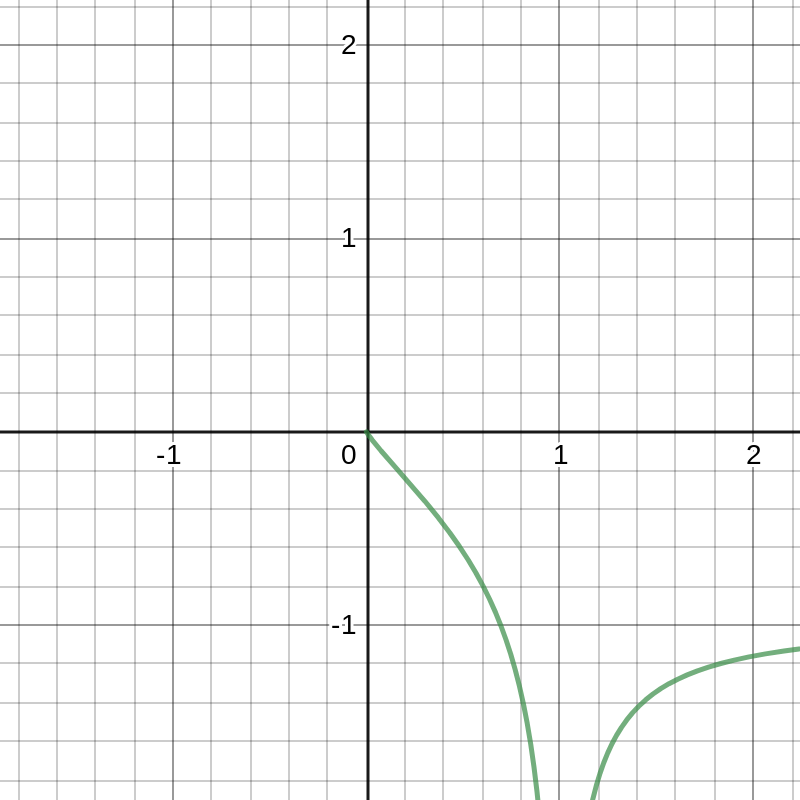}}
	\hspace{1cm}
	\subfloat[ velocity, n=2]{\label{figur:3}\includegraphics[width=50mm]{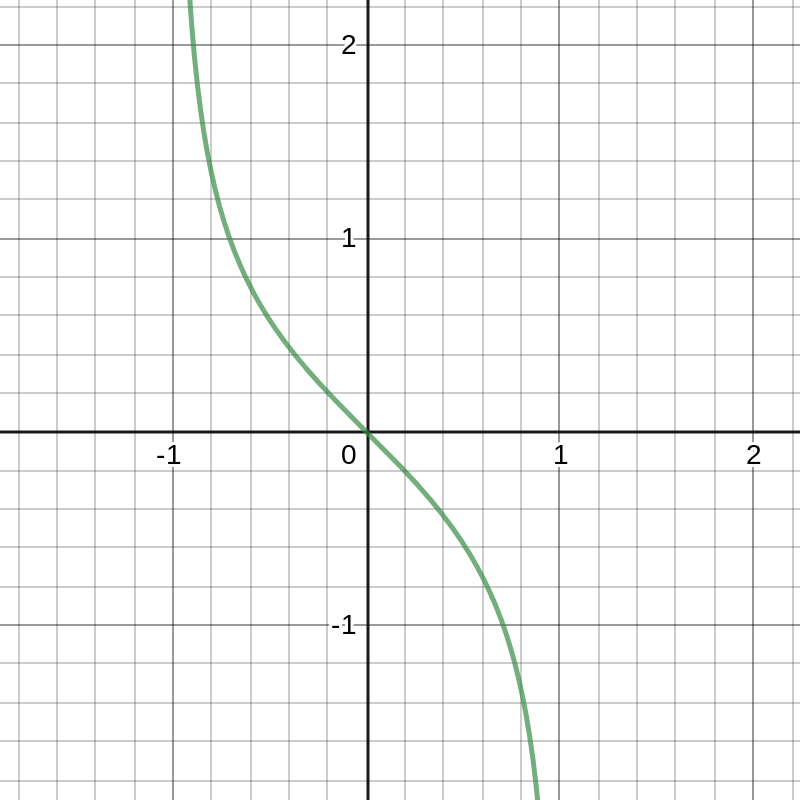}}
	\hspace{1cm}
	\subfloat[ velocity, n=2.1]{\label{figur:3}\includegraphics[width=50mm]{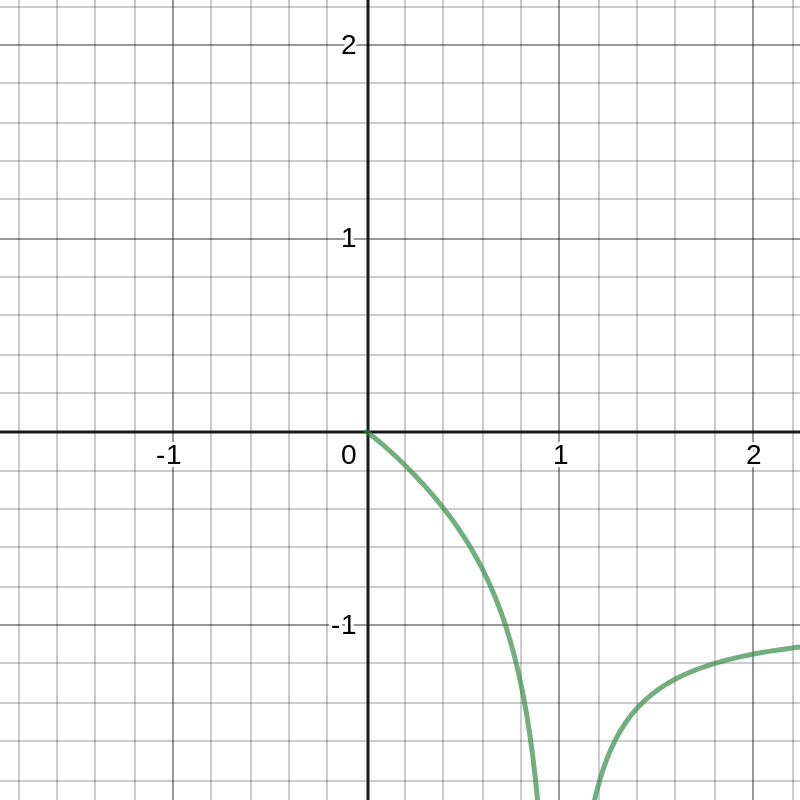}}
	
	\label{figur}\caption{the velocity of angular movement on $\gamma'_{\mathcal{T}}$ in Torus, it shows the behavior near n=2}	
\end{figure*}

\begin{figure*}[htp]
	\centering
	\subfloat[ acceleration, n=1.9 ]{\label{figur:2}\includegraphics[width=50mm]{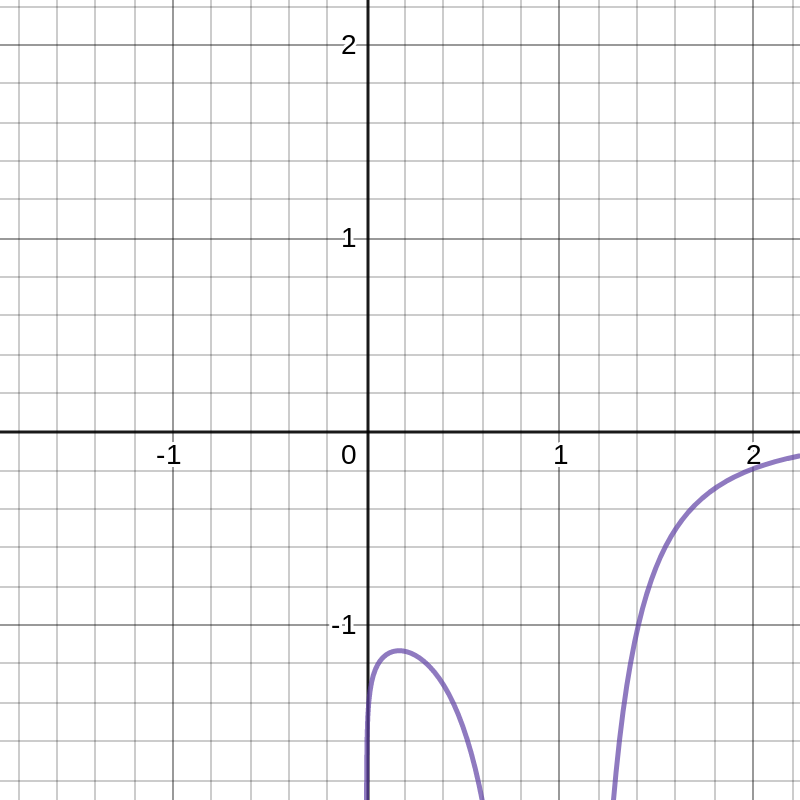}}
	\hspace{1cm}
	\subfloat[ acceleration, n=2]{\label{figur:3}\includegraphics[width=50mm]{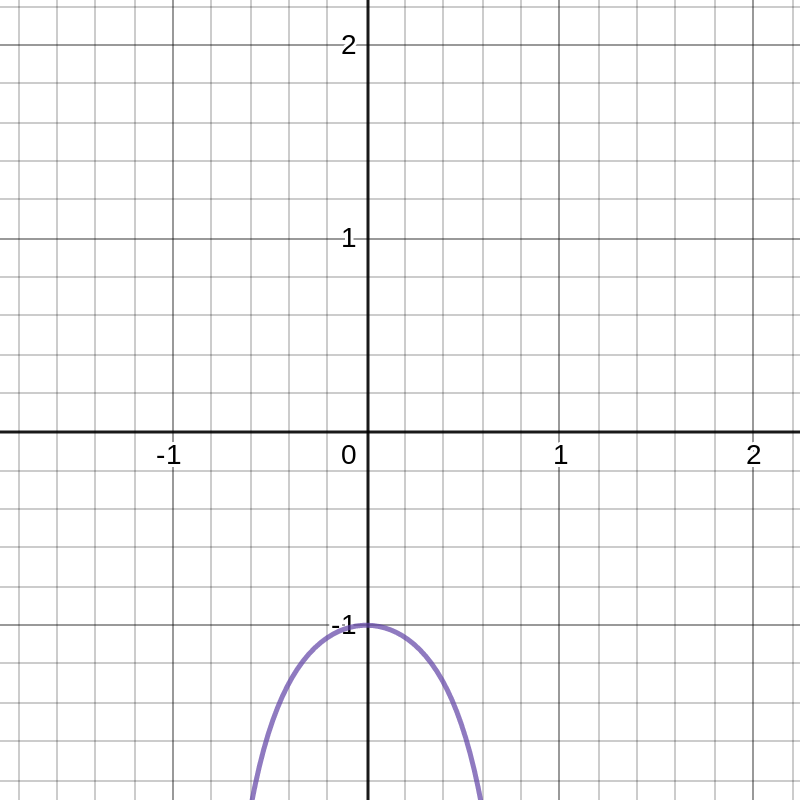}}
	\hspace{1cm}
	\subfloat[ acceleration, n=2.1]{\label{figur:3}\includegraphics[width=50mm]{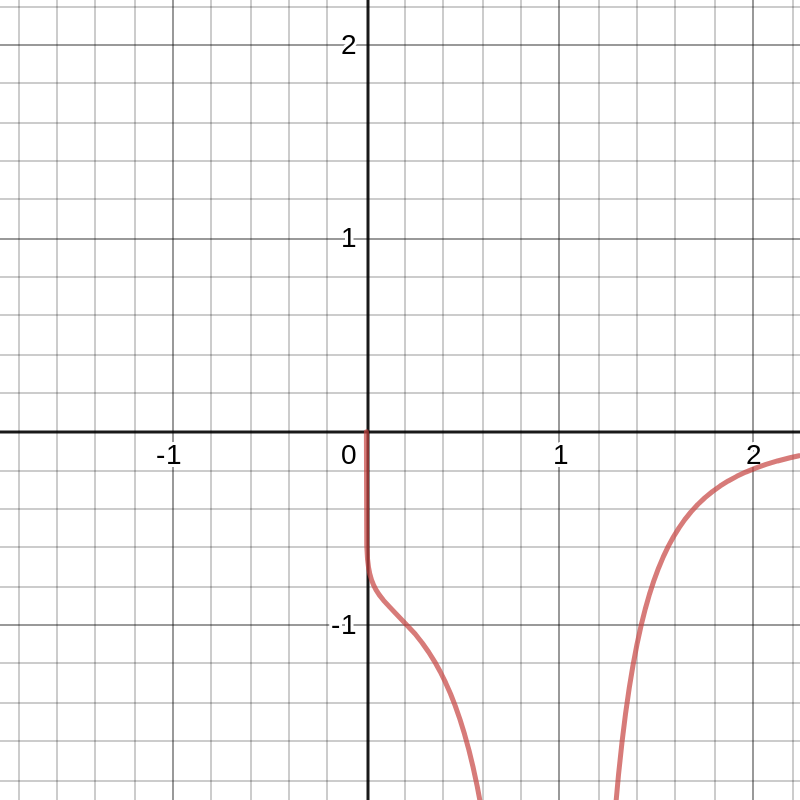}}
	
	\label{figur}\caption{The acceleration of angular movement on $\gamma'_{\mathcal{T}}$ in Torus, it shows the behavior near n=2}	
\end{figure*}

For integrate second derivation, we set new values $w= R+r\cos v$, and then

\begin{equation}\label{eq7}
\dot{w}=-r(\sin v)\dot{v}
\end{equation}

We can divide the first equation of Eq. (33) by $\dot{u}$ and integrate out, then we can get

\begin{equation}\label{eq7}
\int \frac{1}{\dot{u}}\ddot{u}= \int -\frac{2}{w}\dot{w} 
\end{equation}

From Eq. (35) we can get
\begin{equation}\label{eq7}
\dot{u}= \frac{k}{(R + r \cos v)^2} \qed
\end{equation}

Next with similar trick, we can multiply by $\dot{v}$ on the second equation of Eq. (33), we can get Eq. (37)

\begin{equation}\label{eq7}
\int \ddot{v} \dot{v}=\frac{k^2}{a^2} \int \frac{1}{w^3} \dot{w}
\end{equation}

Furthermore, when we integrate out this equation, we can finally get Eq. (38)

\begin{equation}\label{eq7}
\dot{v} = \pm \sqrt{- \frac{k^2}{r^2(R + r\cos u)^2} + l} \qed
\end{equation}

Here $k$ and $l$ are constants of integration. And we want to integrate these equations to get u as a function of v or vice versa. However, as we discuss in \textbf{Section. 2.3} most $(k,l)$ pairs do not yield geodesics.
So when the ratio between $\dot{u}/\dot{v}$ gives integer value, the geodesic of a torus is closed.

\subsection{2.5 Velocity and Acceleration on Torus}

After mapping the straight line on Cartesian space to Torus, the next step is finding a mapped curve for arbitrary line $l'_{\mathcal{U}}$.
In this case, we can also use the same mapping $\mathcal{M}$, like in \textbf{Section. 2.3}.

\begin{align}
\mathcal{M}  &\coloneqq l'_{\mathcal{U}}(a(t),b(t),n) \mapsto \gamma'_{\mathcal{T}}(R(t), r(t), n)
\end{align}

Moreover, one thing we can take a note here is, the mapping for parametrized line $l'_{\mathcal{U}}$, derives the vibration of the Torus radius $r$ and $R$ with time.

We can adjust mapping $\mathcal{M}$ on $\mathcal{F}$, and get the set of projected Fermat's line on Torus.

\begin{gather}
\mathcal{F} \xmapsto{^\mathcal{M}} \mathcal{F}'
\end{gather}

Therefore by \textbf{Theorem. 4}, now we need to find the cross point between set $\mathcal{F}'$ and $\gamma'_\mathcal{T}$.

To find the cross point, we need to calculated moving distance at the point on torus; thus, we need the concept of velocity and acceleration on the point.
In both set $\mathcal{F}$ and $\mathcal{F}'$, the value of the velocity can transform like Eq. (41)

\begin{equation}\label{eq7}
\mathcal{F} \mapsto dy/dx=\frac{dy/dt}{dx/dt}=\frac{dv/dt}{du/dt} \mapsto \mathcal{F}'
\end{equation}

So if we can calculate the derivative of $y$ with $x$ on the Cartesian coordinate, it implies we can derive the velocity at the Torus coordinate.

First, let us expand the equation at \textbf{Theorem. 2} on $\mathbb{R}$ field. Then we can get Eq. (42)

\begin{equation}\label{eq7}
x^n+y^n=1
\end{equation}

We can rewrite $y$ as a function of $x$, in the explicit form Eq. (43), 

\begin{equation}\label{eq7}
y=(1-x^n)^{1/n}
\end{equation}

Also, in the torus, the velocity of $u$ and $v$ can relatively rescale, so we can assume $\dot{u}$=1 without loosing generality. So we can get,

\begin{equation}\label{eq7}
y=\int_{0}^{t} v(t) dt
\end{equation}

\begin{equation}\label{eq7}
vel(t)=\frac{dy}{dx}=-\left(\frac{x}{(1-x^{n})^{1/n}}\right)^{n-1}
\end{equation}

By using Eq. (45), in Fig. 7, you can see the result of $vel(t)$, near n=2.
In the graph, one can see the slope of the $v$ becomes -1, at n=2.

In the same way, we can define acceleration, a(t), like Eq. (46).

\begin{equation}\label{eq7}
a(t)=\frac{d^2y}{dx^2}=-\frac{(n-1)x^{n-2}}{(1-x^{n})^{2-\frac{1}{n}}} 
\end{equation}

In Fig. 8, we can see the result of $a(t)$, near n=2 point. And interestingly,
 near x=0, n=1.9 the acceleration is $-\infty$ but at n=2, x=0, it gives value of -1. Moreover, at n=2.1, the acceleration goes to 0. However, when $x$ becomes slightly larger than 0, the effect of $n$ is negligible.

This is a significant aspect of the acceleration, which could be highly notable evidence for supporting the argument of \textbf{Theorem. 1}.

\section{3. Conclusion and Outlook}

This research suggests alternative form of Fermat's theorem, and shows the acceleration phase transition is occur near n=2, x=0, by using the geometrical mapping which conserves the periodicity of the line element. However, further investigation is needed for the specific relation between this transition and the solution for the Diophantine equation. For example, by using the torus radius vibration, which was shortly discussed in \textbf{Section 2.5}, we can calculate the remainder of distance difference between geodesics and arbitrary parametrized line as a modular of torus inner circumference, 2$\pi$r \qed

\section{4. Appendix}

\subsection{A.1 $b/a$=1, 5 case, Linear mapping on Torus}
In the same process, when $b/q$=1, 5 cases, we can see the mapping result in Fig. 9 and Fig. 10. In the same way every line on a positive integer plane can be a map on the torus with the different $b/a$ values.

\begin{figure}[H]
	\centering
	\includegraphics[angle=0, width=6.5cm, height=3.6cm]{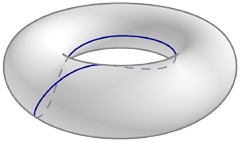}
	\caption{\label{fg3}  The line mapping on Torus. $b/a$=1 case.}
\end{figure}

\begin{figure}[H]
	\centering
	\includegraphics[angle=0, width=6.5cm, height=3.6cm]{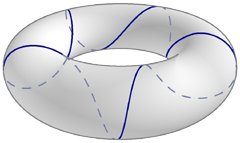}
	\caption{\label{fg3} The line mapping on Torus. $b/a$=5 case.}
\end{figure}

\subsection{A.2 Velocity and Acceleration, at near $n$=3}

In Fig. 11, and Fig. 12, we can observe the behavior of acceleration when $n\neq$ 2.
For example, in Fig. 12, one can see the acceleration result of when n is near 3. Moreover, compare to Fig. 8, there is no phase transition is occurring near x=0 by changing n values.

 \begin{figure*}[htp]
	\centering
	\subfloat[ velocity, n=2.9  ]{\label{figur:1}\includegraphics[width=50mm]{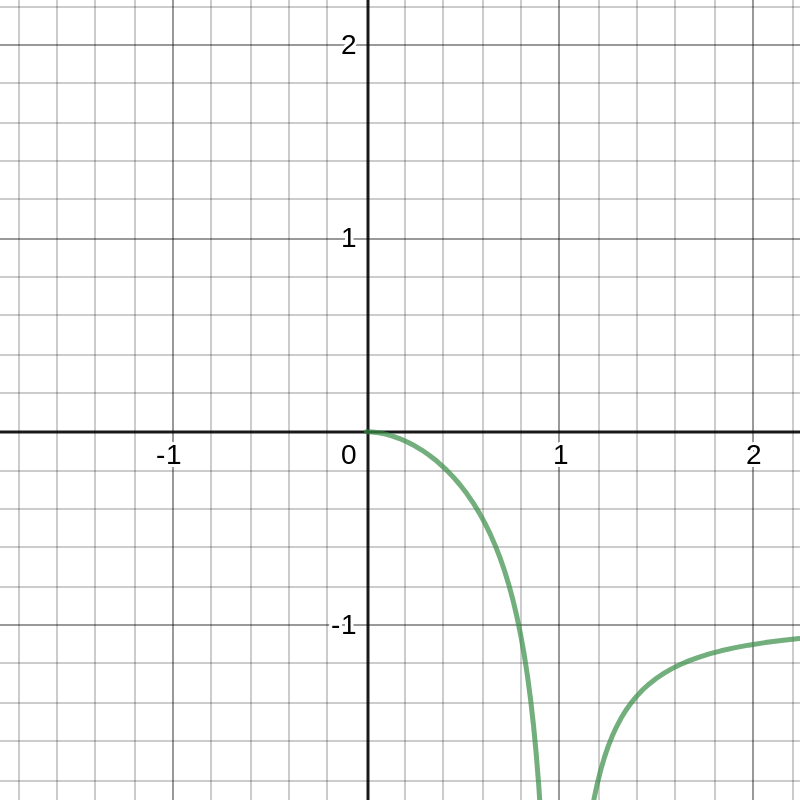}}	
	\hspace{1cm}
	\subfloat[ velocity, n=3 ]{\label{figur:2}\includegraphics[width=50mm]{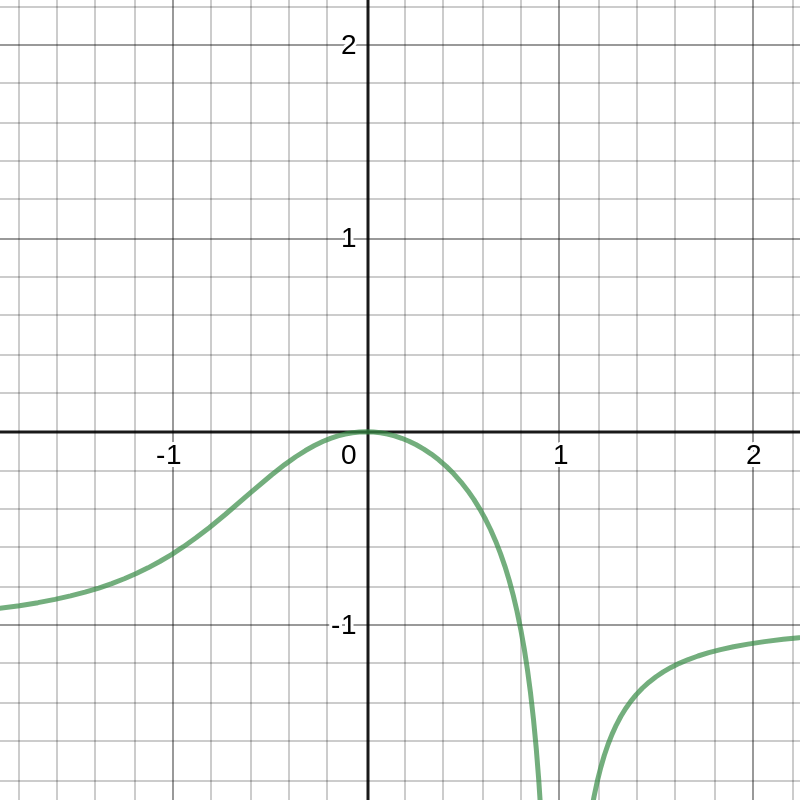}}
	\hspace{1cm}
	\subfloat[ velocity, n=3.1]{\label{figur:3}\includegraphics[width=50mm]{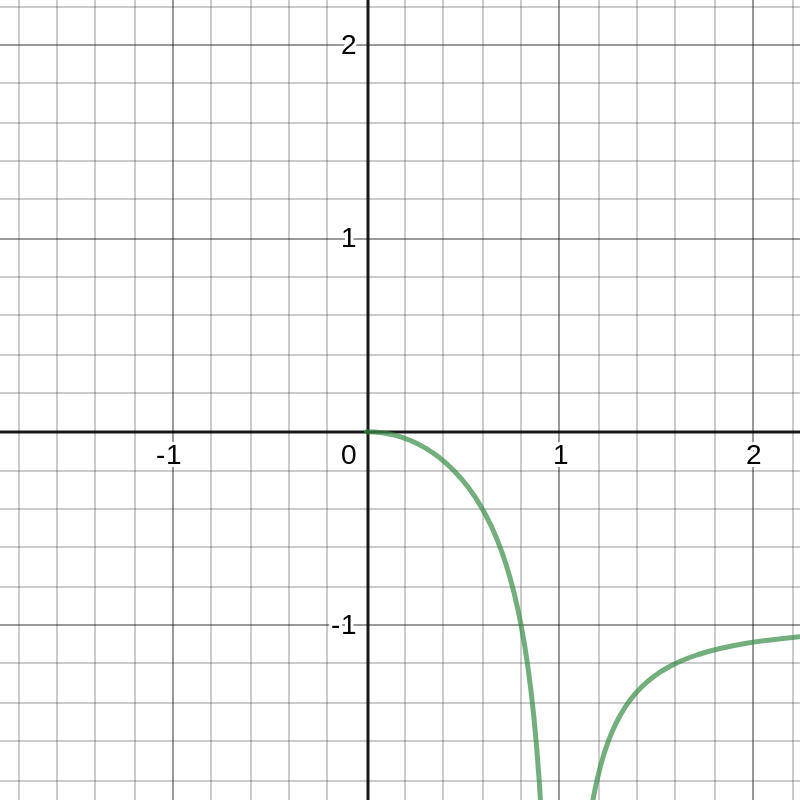}}

	\label{figur}\caption{The velocity of angular movement on $\gamma'_{\mathcal{T}}$ in torus, it shows the behavior near n=3}	
\end{figure*}

\begin{figure*}[htb]
	\centering
	\subfloat[acceleration, n=2.9 ]{\label{figur:1}\includegraphics[width=50mm]{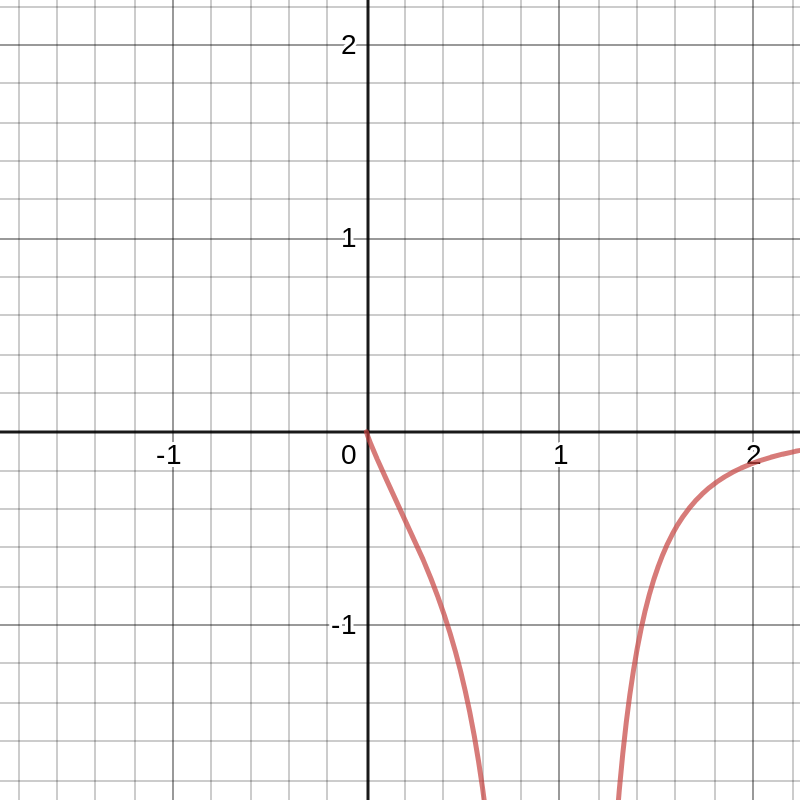}}
	\hspace{1cm}
	\subfloat[acceleration, n=3 ]{\label{figur:2}\includegraphics[width=50mm]{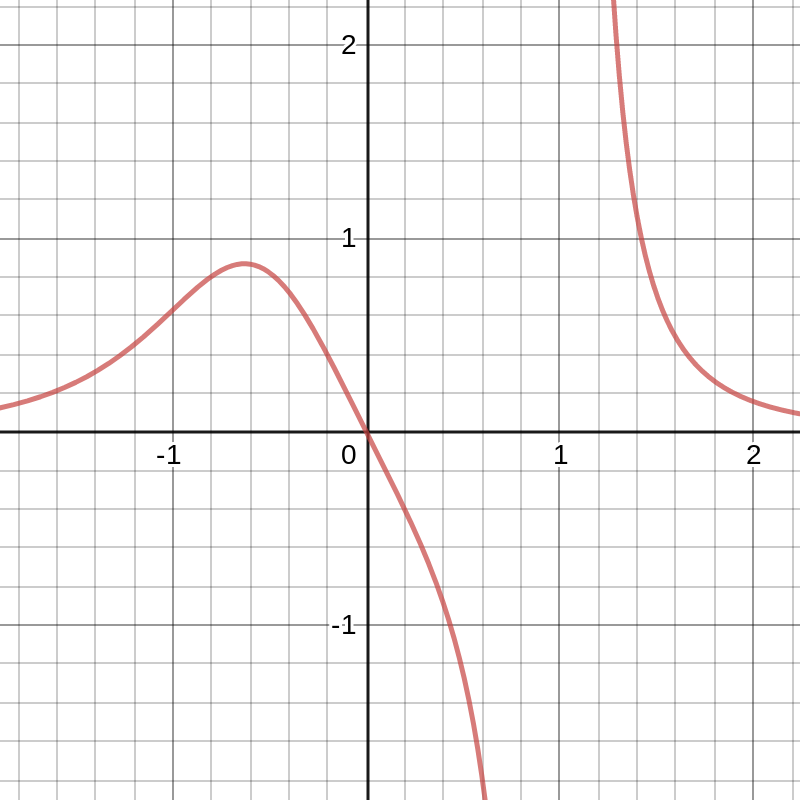}}
	\hspace{1cm}
	\subfloat[acceleration,n=3.1]{\label{figur:3}\includegraphics[width=50mm]{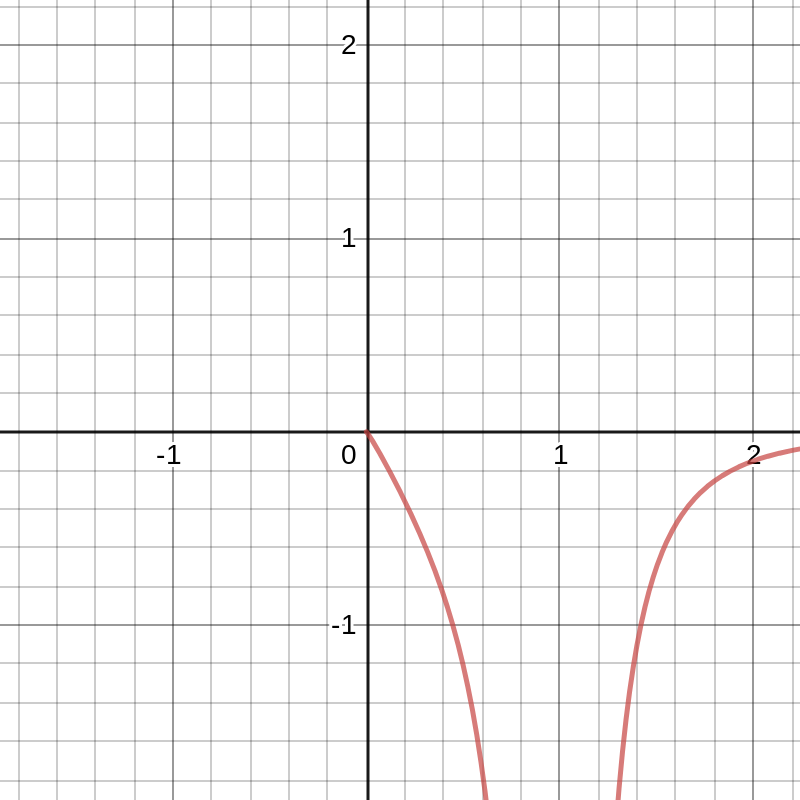}}
	
	\label{figur}\caption{The acceleration of angular movement on $\gamma'_{\mathcal{T}}$ in torus, it shows the behavior near n=3}	
\end{figure*}

\end{document}